\documentclass{article}
\usepackage[utf8]{inputenc}
\usepackage{amsmath,amsthm,amssymb,amsfonts}
\usepackage[usenames]{color}

\usepackage{graphicx}
\usepackage{todonotes}
\usepackage{tikz} 
\usepackage{comment}

\newtheorem{thm}{Theorem}

\def\1{{\bf 1}}

\title{An elementary proof of existence and uniqueness of stationary distributions for irreducible Markov chains
}
\begin{document}

\maketitle

Rinaldo B. Schinazi
\footnote{Department of Mathematics, University of Colorado, Colorado Springs, CO 80933-7150, USA.
E-mail: rinaldo.schinazi@uccs.edu}

\begin{abstract}
Consider an $n\times n$ matrix $P$ with the following properties. All entries in $P$ are positive or $0$, the sum of each row is 1 and for all $i$ and $j$ in $\{1,\dots,n\}$ there exists a natural number $k$ such that the $(i,j)$ entry of the matrix $P^k$ is strictly positive. Then, there exists a unique row vector $v$ with only strictly positive entries, whose sum of entries is 1 and such that $vP=P$. We present a proof of this well-known result that uses only basic algebra and the Bolzano-Weierstrass Theorem.  
\end{abstract}

\section{Introduction}

The existence and uniqueness of stationary distributions for irreducible Markov chains
is a fundamental result in any course on stochastic processes. Its proof is often intertwined with a convergence result, see for instance Bhattharya and Waymire (1990), Feller (1968), Hoel Port and Stone (1972) and Karlin and Taylor (1975). More precisely, the existence of the stationary distribution comes from the convergence of $P^n$ as $n\to\infty$ where $P$ is the transition probability matrix of the Markov chain. While this classical treatment gives more results (convergence in particular) it also requires more hypotheses (aperiodicity of the chain).  More importantly it requires from the reader more technical ability and some probability theory familiarity. It also obscures a linear algebra  result (i.e. a light version of Perron-Frobenius Theorem) which we find interesting in its own right.

In this note we follow the approach of Levin, Peres and Wilmer (2009) to show existence and uniqueness. We use only basic Linear Algebra and the Bolzano-Weierstrass Theorem. The admittedly modest contribution of this note is to streamline the proof and spell out all the details.

\section{The result}

Let $n$ be a natural number and $P$ be an $n\times n$ matrix with the following properties.
\begin{itemize}
    \item (1) $p(i,j)\geq 0$ for all $i$ and $j$ in $\{1,\dots,n\}.$
    \item (2) For all $1\leq i\leq n$, $\sum_{j=1}^n p(i,j)=1$.
    \item (3) For all $i$ and $j$ in $\{1,\dots,n\}$ there exists a natural number $k$ such that the $(i,j)$ entry of the matrix $P^k$ (i.e. the $k$-th power of $P$) is strictly positive.
\end{itemize}

A matrix $P$ with the above properties is a so-called transition probability matrix of an irreducible Markov chain on the finite set $\{1,\dots,n\}$.

\begin{thm}
Consider a matrix $P$ with the three properties listed above. Then, there exists a unique row vector $v$ with only strictly positive entries, whose sum of entries is 1 and such that $vP=P$.
\end{thm}

A row vector $v$ with positive or $0$ entries and whose sum of entries is 1 is a so-called probability distribution. Moreover, if $vP=P$ then $v$ is said to be stationary for the Markov chain associated with matrix $P$. Hence, the Theorem can be rephrased as: An irreducible Markov chain has a unique stationary distribution.

\section{Proof}

{\bf Step 1.} Let $I$ be the $n\times n$ identity matrix. A vector belongs to the kernel of $P-I$
if and only if all its entries are equal. In particular, the dimension of the kernel of the matrix $P-I$ is 1. 

Observe that because of hypothesis (2) any vector whose entries are all equal belongs to the kernel of $P-I$. We now show that these are the only vectors that belong to the kernel.

Let $w$ be a  (column) vector in the kernel of $P-I$. Since $w$ has finitely many entries there exists a natural number $i_0$ in $\{1,\dots,n\}$ such that 
$$\max_{1\leq j\leq n}w(j)=w(i_0).$$

Let $j$ be in $\{1,\dots,n\}$, by hypothesis (3) there exists a natural number $k$ such that the entry  $(i_0,j)$ of the matrix $P^k$ is strictly positive. Hence, (by the definition of matrix multiplication) there is a sequence $i_0,i_1,\dots i_k$ in $\{1,\dots,n\}$ such that $i_k=j$ and
$$p(i_0,i_1)p(i_1,i_2)\dots p(i_{k-1},i_k)>0.$$
By contradiction, assume now that $w(i_1)<w(i_0)$. Since $w=Pw$,
\begin{align*}
w(i_0)=&\sum_{\ell=1}^np(i_0,\ell)w(\ell)\\
   =&p(i_0,i_1)w(i_1)+\sum_{\ell\not =i_1}p(i_0,\ell)w(\ell)\\
    <&p(i_0,i_1)w(i_0)+\sum_{\ell\not =i_1}p(i_0,\ell)w(\ell)\\
    \leq &w(i_0)\left(p(i_0,i_1)+\sum_{\ell\not =i_1}p(i_0,\ell)\right)\\
    =&w(i_0)
\end{align*}
Thus, we get the contradiction $w(i_0)<w(i_0)$. Hence, $w(i_1)=w(i_0)$. By repeating the argument above with $i_1$ and $i_2$ (instead of $i_0$ and $i_1$) we show that $w(i_2)=w(i_1)$. And so on, by successive repetitions of the same argument we show that $w(i_k)=w(i_0)$. But $i_k=j$ and $j$ is arbitrary. Hence, all entries of $w$ are equal. Thus, the only vectors in the kernel of $(P-I)$ are the ones whose entries are all equal. This completes Step 1.

\medskip

{\bf Step 2.} The vector space of all row vectors $v$ such that $vP=v$ is one-dimensional.

By Step 1 we know that the matrix $P-I$ has rank $n-1$. Hence, its transpose $(P-I)'$ has also rank $n-1$. The kernel of $(P-I)'$ is therefore one dimensional. Note now that a row vector $v$ satisfies the equation $vP=v$ if and only if its transpose $v'$ satisfies $P'v'=v'$. That is, there is a one-to-one correspondence between the kernel of $(P-I)'$ and the solutions of the equation $vP=v$. This completes the proof of Step 2.

\medskip

{\bf Step 3.} There exists a unique row vector $\pi$ such that $\pi P=\pi$ with the following properties. All entries are larger than or equal to 0 and the entries sum is 1.

Let $u$ be a row vector whose $n$ entries are all equal to $1/n$. We define a sequence $(v_k)$ of row vectors by
$$v_k=\frac{1}{k}\left(u+uP+\dots uP^{k-1}\right).$$
An easy induction argument shows that for all $j\geq 1$ the matrix $P^j$ has only positive entries and the sum of each row is 1. Thus, for every $j$ the row vector $uP^j$ has only positive entries and the sum of its entries is 1. Therefore, the same is true for $v_k$ for all $k\geq 1$. In particular,
$$||v_k||=\max_{1\leq i\leq n}|v_k(i)|\leq 1.$$
Hence, $(v_k)$ is a bounded sequence in ${\mathbb R}^n$. By Bolzano-Weierstrass Theorem there exists a subsequence $(v_{k_j})$ that converges to some row vector $\pi$. Note that for $1\leq i\leq n$,
$$ |v_kP(i)-v_k(i)|=\frac{1}{k}|uP^k(i)-u(i)|\leq \frac{2}{k}.$$
Hence,
$$ ||v_kP-v_k||=\max_{1\leq i\leq n}|v_kP(i)-v_k(i)|\leq \frac{2}{k}.$$
Therefore, the sequence of row vectors $(v_kP-v_k)$ converges to the zero row vector. The same is true for its subsequence $(v_{k_j}P-v_{k_j})$. Since  $(v_{k_j})$ converges to $\pi$ so does $(v_{k_j}P)$. But (by operations on limits) the sequence $(v_{k_j}P)$ also converges to $\pi P$. Thus, $\pi P=\pi.$ 

Moreover, $\pi$ is the limit of $(v_{k_j})$. Each $v_{k_j}$ has only positive or 0 components and its sum of components is 1. By operations on limits,  the same is true for $\pi$. Since the solutions of the equation $vP=v$ are all in a one-dimensional vector space (Step 2), $\pi$ is the unique probability distribution in this space.
This completes Step 3.

\medskip

{\bf Step 4.} All the entries of the stationary distribution $\pi$ are strictly positive.

First observe that $\pi P=\pi$. Hence, $\pi P^2=\pi P=\pi$. More generally, for every natural number $k\geq 1$, $\pi P^k=\pi$. 

Assume by contradiction that there is $i$ in $\{1,\dots,n\}$ such that $\pi(i)=0.$
Let $j$ be in $\{1,\dots,n\}$. By Hypothesis (3) there is a natural number $k\geq 1$ such that $p_k(j,i)>0$, where $p_k(j,i)$ is the $(j,i$) entry of matrix $P^k$. Since $\pi P^k=\pi$,
\begin{align*}
    \pi(i)=&\sum_{\ell=1}^n \pi(\ell)p_k(\ell,i)\\
    \geq &\pi(j)p_k(j,i)\\
\end{align*}
But $\pi(i)=0$ and $p_k(j,i)>0$ therefore $\pi(j)=0$. Since this is true for all $j$ in $\{1,\dots,n\}$, $\pi$ is the zero row vector. By Step 3 we know this is not so. We have a contradiction. Hence, $\pi(i)>0$ for all $i$ in $\{1,\dots,n\}$. This completes Step 4 and the proof of the Theorem.

\bibliographystyle{amsplain}

\end{document}